\documentclass{amsart}
\usepackage[latin1]{inputenc}
\usepackage[T1]{fontenc}
\usepackage[english]{babel}
\usepackage{subfigure}
\usepackage{graphicx}

\newtheorem{assumption}{Assumption}
\newtheorem{theorem}{Theorem}
\newtheorem{remark}{Remark}

\begin{document}

\title[A remark on the boundedness and convergence properties]{A remark on the boundedness and 
convergence properties of smooth sliding mode controllers}
\author{Wallace M. Bessa}
\address{CEFET/RJ, Federal Center for Technological Education, Av.\ Maracanã 229, DEPES, CEP 
20271-110, Rio de Janeiro, Brazil}
\email{wmbessa@cefet-rj.br, wmbessa@ams.org}
\subjclass[2000]{Primary 93C10, 93C40; Secondary 93C15, 34H05}
\keywords{Convergence analysis; Lyapunov methods; Nonlinear control; Sliding mode.}
\begin{abstract}
Conventional sliding mode controllers are based on the assumption of switching control but a 
well-known drawback of this approach is the chattering phenomenon. To overcome the undesirable 
chattering effects, the discontinuity in the control law can be smoothed out in a thin boundary 
layer neighboring the switching surface. In this work, rigorous proofs of the boundedeness and 
convergence properties of smooth sliding mode controllers are presented. This result corrects 
flawed conclusions previously reached in the literature. 
\end{abstract}

\maketitle

\section{Introduction}

Sliding mode control theory was conceived and developed in the former Soviet Union by Emelyanov
\cite{emelyanov}, Filippov \cite{filippov}, Itkis \cite{itkis}, Utkin \cite{utkin} and others. 
But a known drawback of conventional sliding mode controllers is the chattering phenomenon due to
the discontinuous term in the control law. In order to avoid the undesired effects of the control 
chattering, Slotine \cite{slotine} proposed the adoption of a thin boundary layer neighboring the 
switching surface, by replacing the sign function by a saturation function. This substitution can 
minimize or, when desired, even completely eliminate chattering, but turns {\em perfect tracking} 
into a {\em tracking with guaranteed precision} problem, which actually means that a steady-state 
error will always remain. 

This paper presents a convergence analysis of smooth sliding mode controllers. The finite-time
convergence of the tracking error vector to the boundary layer is handled using Lyapunov's direct 
method. It is also analytically proven that, once in boundary layer, the error vector exponentially 
converges to a bounded region. This result corrects a minor flaw in Slotine's work, by showing that 
the tracking error bounds are different from the bounds provided in \cite{slotine}.

\section{Problem statement and controller design}

Consider a class of $n^\mathrm{th}$-order nonlinear system:

\begin{equation}
x^{(n)}=f(\mathbf{x})+b(\mathbf{x})u\\
\label{eq:system}
\end{equation}

\noindent
where $u$ is the control input, the scalar variable $x$ is the output of interest, $x^{(n)}$ is the 
$n^\mathrm{th}$ derivative of $x$ with respect to time $t\in[0,\infty)$, $\mathbf{x}=[x,\dot{x},
\ldots,x^{(n-1)}]$ is the system state vector, and $f,b:\mathbb{R}^n\rightarrow\mathbb{R}$ are 
both nonlinear functions.

In respect of the dynamic system presented in Eq.~(\ref{eq:system}), the following assumptions will 
be made:

\begin{assumption}
The function $f$ is unknown but bounded by a known function of $\mathbf{x}$, i.e., $|\hat{f}
(\mathbf{x})-f(\mathbf{x})|\le F(\mathbf{x})$ where $\hat{f}$ is an estimate of $f$.
\label{as:fbounds}
\end{assumption}
\begin{assumption}
The input gain $b(\mathbf{x})$ is unknown but positive and bounded, i.e., $0<b_{\mathrm{min}}
\le b(\mathbf{x})\le b_{\mathrm{max}}$.
\label{as:bbounds}
\end{assumption}

The proposed control problem is to ensure that, even in the presence of parametric uncertainties and
unmodeled dynamics, the state vector $\mathbf{x}$ will follow a desired trajectory $\mathbf{x}_d=
[x_d,\dot{x}_d,\ldots,x^{(n-1)}_d]$ in the state space.

Regarding the development of the control law the following assumptions should also be made:

\begin{assumption}
The state vector $\mathbf{x}$ is available.
\label{as:stat}
\end{assumption}
\begin{assumption}
The desired trajectory $\mathbf{x}_d$ is once differentiable in time. Furthermore, every element
of vector $\mathbf{x}_d$, as well as $x^{(n)}_d$, is available and with known bounds.
\label{as:traj}
\end{assumption}

Now, let $\tilde{x}=x-x_d$ be defined as the tracking error in the variable $x$, and 

\begin{displaymath}
\mathbf{\tilde{x}}=\mathbf{x}-\mathbf{x}_d=[\tilde{x},\dot{\tilde{x}},\ldots,\tilde{x}^{(n-1)}]
\end{displaymath}

\noindent 
as the tracking error vector. 

Consider a sliding surface $S$ defined in the state space by the equation $s(\mathbf{\tilde{x}})=0$, with 
the function $s:\mathbb{R}^n\rightarrow\mathbb{R}$ satisfying

\begin{displaymath}
\displaystyle s(\mathbf{\tilde{x}})=\left(\frac{d}{dt}+\lambda\right)^{n-1}\tilde{x}
\end{displaymath}

\noindent 
or conveniently rewritten as

\begin{equation}
s(\mathbf{\tilde{x}})=\mathbf{c^\mathrm{T}\tilde{x}}
\label{eq:surf}
\end{equation}

\noindent
where $\mathbf{c}=[c_{n-1}\lambda^{n-1},\ldots,c_1\lambda,c_0]$ and $c_i$ states for 
binomial coefficients, i.e.,

\begin{equation}
c_i=\binom{n-1}{i}=\frac{(n-1)!}{(n-i-1)!\:i!}\:,\quad i=0,1,\ldots,n-1 
\label{eq:binom}
\end{equation}

\noindent
which makes $c_{n-1}\lambda^{n-1}+\cdots+c_1\lambda+c_0$ a Hurwitz polynomial. 

From Eq.~(\ref{eq:binom}), it can be easily verified that $c_0=1$, for $\forall n\ge1$. Thus, for 
notational convenience, the time derivative of $s$ will be written in the following form:

\begin{equation}
\dot{s}=\mathbf{c^\mathrm{T}\dot{\tilde{x}}}
=\tilde{x}^{(n)}+\mathbf{\bar{c}^\mathrm{T}\tilde{x}}
\label{eq:sd}
\end{equation}

\noindent
where $\mathbf{\bar{c}}=[0,c_{n-1}\lambda^{n-1},\ldots,c_1\lambda]$.

Now, let the problem of controlling the uncertain nonlinear system (\ref{eq:system}) be treated 
via the classical sliding mode approach, defining a control law composed by an equivalent control 
$\hat{u}=\hat{b}^{-1}(-\hat{f}+x^{(n)}_d-\mathbf{\bar{c}^\mathrm{T}\tilde{x}})$ and a 
discontinuous term $-K\mbox{sgn}(s)$:

\begin{equation}
u=\hat{b}^{-1}\left(-\hat{f}+x^{(n)}_d-\mathbf{\bar{c}^\mathrm{T}\tilde{x}}
\right)-K\mbox{sgn}(s)
\label{eq:usgn}
\end{equation}

\noindent
where $\hat{b}=\sqrt{b_\mathrm{max}b_\mathrm{min}}$ is an estimate of $b$, $K$ is a positive gain 
and sgn($\cdot$) is defined as 

\begin{displaymath}
\mbox{sgn}(s) = \left\{\begin{array}{rc}
-1&\mbox{if}\quad s<0 \\
0&\mbox{if}\quad s=0 \\
1&\mbox{if}\quad s>0
\end{array}\right.
\end{displaymath}

Based on Assumptions~\ref{as:fbounds} and \ref{as:bbounds} and considering that $\beta^{-1}\le
\hat{b}/b\le\beta$, where $\beta=\sqrt{b_\mathrm{max}/b_\mathrm{min}}$, the gain $K$ should be 
chosen according to

\begin{equation}
K\ge\beta\hat{b}^{-1}(\eta+F)+(\beta-1)\big|\hat{b}^{-1}(-\hat{f}+x^{(n)}_d
-\mathbf{\bar{c}^\mathrm{T}\tilde{x}})\big|
\label{eq:gain}
\end{equation}

\noindent
where $\eta$ is a strictly positive constant related to the reaching time. 

Therefore, it can be easily verified that (\ref{eq:usgn}) is sufficient to impose the sliding condition 

\begin{displaymath}
\displaystyle\frac{1}{2}\frac{d}{dt}s^2\le-\eta|s|
\end{displaymath}

\noindent
which, in fact, ensures the finite-time convergence of the tracking error vector to the sliding 
surface $S$ and, consequently, its exponential stability.

However, the presence of a discontinuous term in the control law leads to the well known chattering 
effect. To avoid these undesirable high-frequency oscillations of the controlled variable, Slotine 
\cite{slotine} proposed the adoption of a a thin boundary layer, $S_\phi$, in the neighborhood of 
the switching surface:

\begin{equation}
S_\phi=\big\{\mathbf{x}\in\mathbb{R}^n\:\big|\:|s(\mathbf{\tilde{x}})|\le\phi\big\}
\label{eq:bound}
\end{equation}

\noindent
where $\phi$ is a strictly positive constant that represents the boundary layer thickness.

The boundary layer is achieved by replacing the sign function by a continuous interpolation inside 
$S_\phi$. It should be emphasized that this smooth approximation, which will be called here 
$\varphi(s,\phi)$, must behave exactly like the sign function outside the boundary layer. There 
are several options to smooth out the ideal relay but the most common choices are the saturation 
function:

\begin{equation}
\mbox{sat}(s/\phi) = \left\{\begin{array}{cc}
\mbox{sgn}(s)&\mbox{if}\quad |s/\phi|\ge1 \\
s/\phi&\mbox{if}\quad |s/\phi|<1 
\end{array}\right.
\label{eq:sat}
\end{equation}

\noindent
and the hyperbolic tangent function $\tanh(s/\phi)$.

In this way, the smooth sliding mode control law can be stated as follows

\begin{equation}
u=\hat{b}^{-1}\left(-\hat{f}+x^{(n)}_d-\mathbf{\bar{c}^\mathrm{T}\tilde{x}}\right)-K\varphi(s,\phi)
\label{eq:usat}
\end{equation}

\section{Convergence analysis}

The attractiveness and invariant properties of the boundary layer are established in the following 
theorem.

\begin{theorem}
\label{th:theo1}
Consider the uncertain nonlinear system (\ref{eq:system}) and Assumptions~\ref{as:fbounds}--\ref{%
as:traj}. Then, the smooth sliding mode controller defined by (\ref{eq:usat}) and (\ref{eq:gain}) 
ensures the finite-time convergence of the tracking error vector to the boundary layer $S_\phi$, 
defined according to (\ref{eq:bound}).
\end{theorem}

{\bf Proof:} 
Let a positive-definite Lyapunov function candidate $V$ be defined as

\begin{displaymath}
\displaystyle V(t)=\frac{1}{2}s^2_\phi
\end{displaymath}

\noindent
where $s_\phi$ is a measure of the distance of the current error to the boundary layer, and can be 
computed as follows

\begin{equation}
s_\phi=s-\phi\:\mbox{sat}(s/\phi)
\label{eq:dist}
\end{equation}

Noting that $s_\phi=0$ in the boundary layer, one has $\dot{V}(t)=0$ inside $S_\phi$. From Eqs.
(\ref{eq:sat}) and (\ref{eq:dist}), it can be easily verified that $\dot{s}_\phi=\dot{s}$ outside 
the boundary layer and, in this case, $\dot{V}$ becomes

\begin{align*}
\dot{V}(t)&=s_\phi\dot{s}_\phi=s_\phi\dot{s}=(x^{(n)}-x^{(n)}_d+\mathbf{\bar{c}^\mathrm{T}
\tilde{x}})s_\phi\\
&=\left(f+bu-x^{(n)}_d+\mathbf{\bar{c}^\mathrm{T}\tilde{x}}\right)s_\phi
\end{align*}

Considering that outside the boundary layer the control law (\ref{eq:usat}) takes the following 
form:

\begin{displaymath}
u=\hat{b}^{-1}\left(-\hat{f}+x^{(n)}_d-\mathbf{\bar{c}^\mathrm{T}\tilde{x}}\right)
-K\mbox{sgn}(s_\phi)
\end{displaymath}

\noindent
and noting that $f=\hat{f}-(\hat{f}-f)$, one has

\begin{align*}
\dot{V}(t)=&\big[f+b\hat{b}^{-1}(-\hat{f}+x^{(n)}_d-\mathbf{\bar{c}^\mathrm{T}\tilde{x}})
-bK\mbox{sgn}(s_\phi)-x^{(n)}_d+\mathbf{\bar{c}^\mathrm{T}\tilde{x}}\big]s_\phi\\
=&-\big[(\hat{f}-f)-b\hat{b}^{-1}(-\hat{f}+x^{(n)}_d-\mathbf{\bar{c}^\mathrm{T}\tilde{x}})
+(-\hat{f}+x^{(n)}_d-\mathbf{\bar{c}^\mathrm{T}\tilde{x}})+bK\mbox{sgn}(s_\phi)\big]s_\phi
\end{align*}

So, considering Assumptions~\ref{as:fbounds} and \ref{as:bbounds} and defining $K$ according to 
(\ref{eq:gain}), $\dot{V}$ becomes:

\begin{displaymath}
\dot{V}(t)\le-\eta|s_\phi|
\end{displaymath}

\noindent
which implies $V(t)\le V(0)$ and that $s_\phi$ is bounded. From the definition of $s_\phi$, it can 
be easily verified that $s$ is bounded. Considering Assumption~\ref{as:traj} and Eq.~(\ref{eq:sd}), 
it can be concluded that $\dot{s}$ is also bounded.

The finite-time convergence of the tracking error vector to the boundary layer can be shown by 
recalling that

\begin{displaymath}
\dot{V}(t)=\frac{1}{2}\frac{d}{dt}s_{\phi}^2=s_\phi\dot{s}_\phi\le-\eta|s_\phi|  
\end{displaymath}

Then, dividing by $|s_\phi|$ and integrating both sides between 0 and $t$ gives

\begin{displaymath}
\int^t_0\frac{s_\phi}{|s_\phi|}\dot{s}_\phi\,d\tau\le-\int^t_0\eta\,d\tau
\end{displaymath}

\begin{displaymath}
|s_\phi(t)|-|s_\phi(0)|\le-\eta\,t 
\end{displaymath}

In this way, considering $t_\mathrm{reach}$ as the time required to hit $S_\phi$ and noting that 
$|s_\phi(t_\mathrm{reach})|=0$, one has

\begin{displaymath}
\displaystyle t_{\mathrm{reach}}\le\frac{|s_\phi(0)|}{\eta}
\end{displaymath}

\noindent
which guarantees the convergence of the tracking error vector to the boundary layer in a time 
interval smaller than $|s_\phi(0)|/\eta$ and completes the proof.
\hfill$\square$\vspace*{10pt}

\begin{figure}[htb]
\centering
\includegraphics[width=0.4\textwidth]{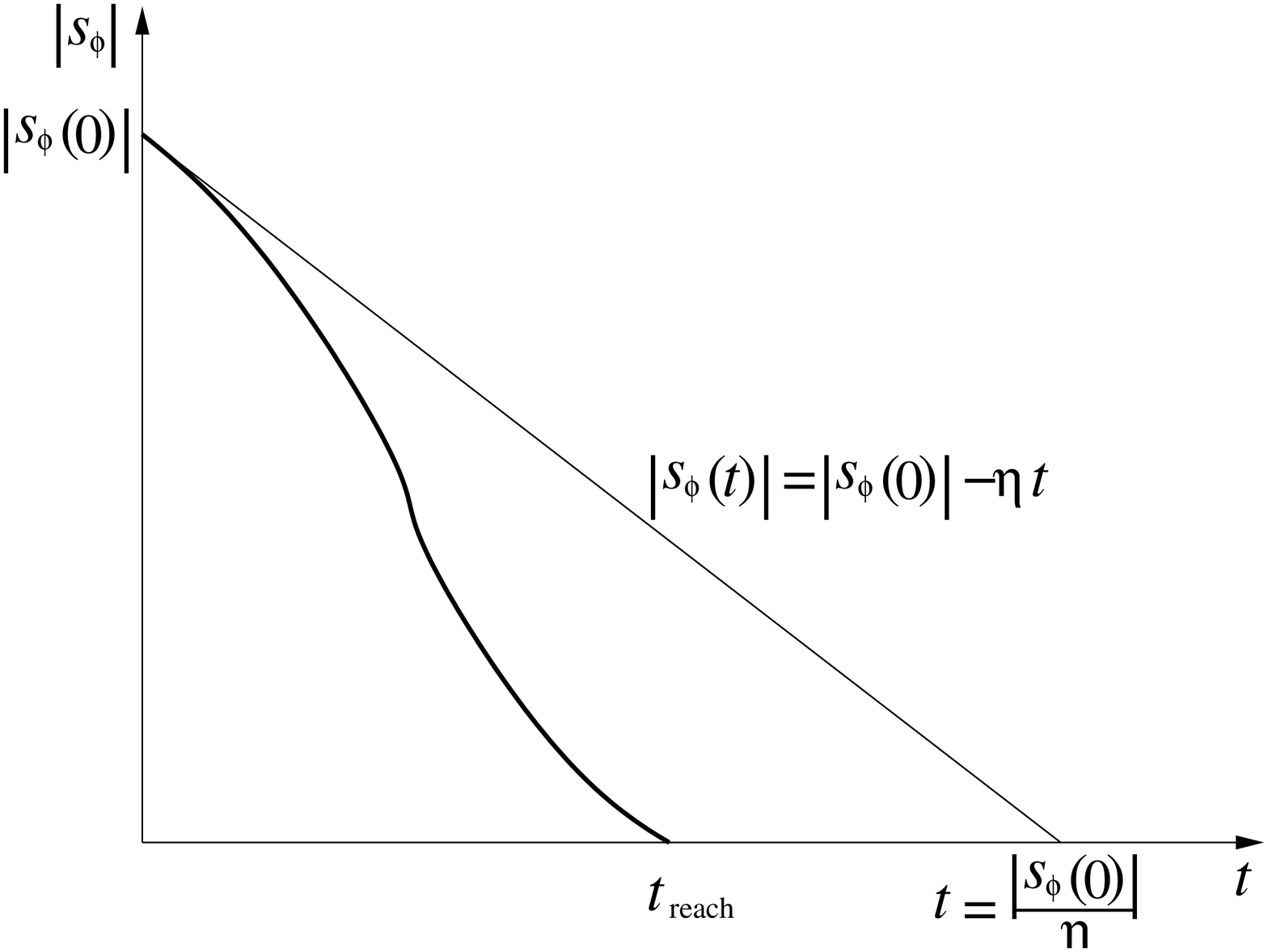} 
\caption{Time evolution of $|s_\phi|$.}
\label{fi:treach}
\end{figure}

Therefore, the value of the positive constant $\eta$ can be properly chosen in order to keep the 
reaching time, $t_\mathrm{reach}$, as short as possible. Figure~\ref{fi:treach} shows that the 
time evolution of $|s_\phi|$ is bounded by the straight line $|s_\phi(t)|=|s_\phi(0)|-\eta\,t$.

Finally, the proof of the boundedness of the tracking error vector relies on Theorem~\ref{th:theo2}.

\begin{theorem}
\label{th:theo2}
Let the boundary layer $S_\phi$ be defined according to (\ref{eq:bound}). Then, once inside 
$S_\phi$, the tracking error vector will exponentially converge to a closed region $\Phi=\{
\mathbf{x}\in\mathbb{R}^n\:|\:|s(\mathbf{\tilde{x}})|\le\phi\mbox{ and }|\tilde{x}^{(i)}|\le\zeta_i
\lambda^{i-n+1}\phi,i=0,1,\ldots,n-1\}$, with $\zeta_i$ defined as

\begin{equation}
\zeta_i = \left\{\begin{array}{cl}
1&\mbox{for}\quad i=0 \\
1+\sum^{i-1}_{j=0}\binom{i}{j}\zeta_j&\mbox{for}\quad i=1,2,\ldots,n-1.
\end{array}\right.
\label{eq:zeta}
\end{equation}

\end{theorem}

{\bf Proof:} 
From the definition of $s$, Eq.~(\ref{eq:surf}), and considering that $|s(\mathbf{x})|\le\phi$ 
may be rewritten as $-\phi\le s(\mathbf{x})\le\phi$, one has 

\begin{equation}
-\phi\le c_0\tilde{x}^{(n-1)}+c_1\lambda\tilde{x}^{(n-2)}+\cdots+c_{n-1}\lambda^{n-1}\tilde{x}\le\phi
\label{eq:sbounds}
\end{equation}

Multiplying (\ref{eq:sbounds}) by $e^{\lambda t}$ yields

\begin{equation}
-\phi e^{\lambda t}\le\frac{d^{n-1}}{dt^{n-1}}(\tilde{x}e^{\lambda t})\le\phi e^{\lambda t} 
\label{eq:dbounds}
\end{equation}

Integrating (\ref{eq:dbounds}) between $0$ and $t$ gives

\begin{equation}
-\frac{\phi}{\lambda} e^{\lambda t}+\frac{\phi}{\lambda}\le\frac{d^{n-2}}{dt^{n-2}}(\tilde{x}
e^{\lambda t})-\left.\frac{d^{n-2}}{dt^{n-2}}(\tilde{x}e^{\lambda t})\right|_{t=0}\le
\frac{\phi}{\lambda} e^{\lambda t}-\frac{\phi}{\lambda} 
\label{eq:int_1a}
\end{equation}

\noindent
or conveniently rewritten as

\begin{multline}
-\frac{\phi}{\lambda} e^{\lambda t}-\left(\left|\frac{d^{n-2}}{dt^{n-2}}(\tilde{x}e^{\lambda t})
\right|_{t=0}+\frac{\phi}{\lambda}\right)\\\le\frac{d^{n-2}}{dt^{n-2}}(\tilde{x}e^{\lambda t})\le\\
\frac{\phi}{\lambda} e^{\lambda t}+\left(\left|\frac{d^{n-2}}{dt^{n-2}}(\tilde{x}e^{\lambda t})
\right|_{t=0}+\frac{\phi}{\lambda}\right)
\label{eq:int_1b}
\end{multline}

The same reasoning can be repeatedly applied until the $(n-1)^\mathrm{th}$ integral of 
(\ref{eq:dbounds}) is reached:

\begin{multline}
-\frac{\phi}{\lambda^{n-1}}e^{\lambda t}-\left(\left|\frac{d^{n-2}}{dt^{n-2}}(\tilde{x}e^{\lambda
t})\right|_{t=0}+\frac{\phi}{\lambda}\right)\frac{t^{n-2}}{(n-2)!}-\cdots\\-\left(\left|\tilde{x}
(0)\right|+\frac{\phi}{\lambda^{n-1}}\right)\le\tilde{x}e^{\lambda t}\le\frac{\phi}{\lambda^{n-1}}
e^{\lambda t}+\\\left(\left|\frac{d^{n-2}}{dt^{n-2}}(\tilde{x}e^{\lambda t})\right|_{t=0}
+\frac{\phi}{\lambda}\right)\frac{t^{n-2}}{(n-2)!}+\cdots+\left(\left|\tilde{x}(0)\right|
+\frac{\phi}{\lambda^{n-1}}\right)
\label{eq:int_n-1}
\end{multline}

Furthermore, dividing (\ref{eq:int_n-1}) by $e^{\lambda t}$, it can be easily verified that, for 
$t\to\infty$,

\begin{equation}
-\frac{\phi}{\lambda^{n-1}}\le\tilde{x}(t)\le\frac{\phi}{\lambda^{n-1}}
\label{eq:txbound}
\end{equation}

Considering the $(n-2)^\mathrm{th}$ integral of (\ref{eq:dbounds})

\begin{multline}
-\frac{\phi}{\lambda^{n-2}}e^{\lambda t}-\left(\left|\frac{d^{n-2}}{dt^{n-2}}(\tilde{x}
e^{\lambda t})\right|_{t=0}+\frac{\phi}{\lambda}\right)\frac{t^{n-3}}{(n-3)!}-\cdots\\
-\left(\left|\dot{\tilde{x}}(0)\right|+\frac{\phi}{\lambda^{n-2}}\right)\le\frac{d}{dt}
(\tilde{x}e^{\lambda t})\le\frac{\phi}{\lambda^{n-2}}e^{\lambda t}+\\\left(\left|
\frac{d^{n-2}}{dt^{n-2}}(\tilde{x}e^{\lambda t})\right|_{t=0}+\frac{\phi}{\lambda}\right)
\frac{t^{n-3}}{(n-3)!}+\cdots+\left(\left|\dot{\tilde{x}}(0)\right|+\frac{\phi}{\lambda^{n-2}}
\right)
\label{eq:int_n-2}
\end{multline}

\noindent
and noting that $d(\tilde{x}e^{\lambda t})/dt=\dot{\tilde{x}}e^{\lambda t}+\tilde{x}\lambda e^{\lambda
t}$, by imposing the bounds (\ref{eq:txbound}) to (\ref{eq:int_n-2}) and dividing again by $e^{\lambda
t}$, it follows that, for $t\to\infty$,
 
\begin{equation}
-2\frac{\phi}{\lambda^{n-2}}\le\dot{\tilde{x}}(t)\le2\frac{\phi}{\lambda^{n-2}}
\label{eq:txdbound}
\end{equation}

Now, applying the bounds (\ref{eq:txbound}) and (\ref{eq:txdbound}) to the $(n-3)^\mathrm{th}$ integral 
of (\ref{eq:dbounds}) and dividing once again by $e^{\lambda t}$, it follows that, 
for $t\to\infty$,

\begin{equation}
-6\frac{\phi}{\lambda^{n-3}}\le\ddot{\tilde{x}}(t)\le6\frac{\phi}{\lambda^{n-3}}
\label{eq:txddbound}
\end{equation}

The same procedure can be successively repeated until the bounds for $\tilde{x}^{(n-1)}$ are 
achieved:

\begin{equation}
-\left(1+\sum^{n-2}_{i=0}\binom{n-1}{i}\zeta_i\right)\phi\le\tilde{x}^{(n-1)}\le
\left(1+\sum^{n-2}_{i=0}\binom{n-1}{i}\zeta_i\right)\phi
\label{eq:txnbound}
\end{equation}

\noindent
where the coefficients $\zeta_i$ ($i=0,1,\ldots,n-2$) are related to the previously obtained bounds 
of each $\tilde{x}^{(i)}$ and can be summarized as in (\ref{eq:zeta}).

In this way, by inspection of the integrals of (\ref{eq:dbounds}), as well as (\ref{eq:txbound}), 
(\ref{eq:txdbound}), (\ref{eq:txddbound}), (\ref{eq:txnbound}) and the other omitted bounds, it follows 
that the tracking error will be confined within the limits $|\tilde{x}^{(i)}|\le\zeta_i\lambda^{i-n+1}
\phi,i=0,1,\ldots,n-1$, where $\zeta_i$ is defined by (\ref{eq:zeta}). 

\begin{figure}[htb]
\centering
\includegraphics[width=0.4\textwidth]{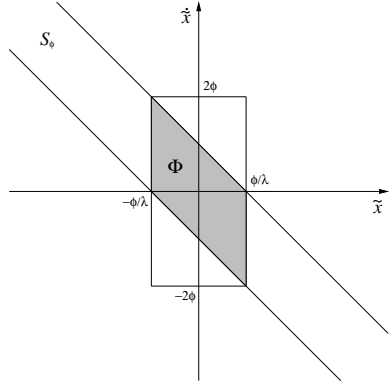} 
\caption{Bounds of $x^{(i)}$ for a $2^\mathrm{nd}$-order system.}
\label{fi:bounds}
\end{figure}

However, the aforementioned bounds define an $n$-dimensional box that is not completely inside the 
boundary layer. Figure.~\ref{fi:bounds} illustrates for a $2^\mathrm{nd}$-order system ($n=2$).

Considering the attractiveness and invariant properties of $S_\phi$ demonstrated in Theorem 
\ref{th:theo1}, the region of convergence can be stated as the intersection of the boundary layer 
and the $n$-dimensional box defined by the preceding bounds. Therefore, it follows that the tracking 
error vector will exponentially converge to a closed region $\Phi=\{\mathbf{x}\in\mathbb{R}^n\:|\:
|s(\mathbf{\tilde{x}})|\le\phi\mbox{ and }|\tilde{x}^{(i)}|\le\zeta_i\lambda^{i-n+1}\phi,i=0,1,
\ldots,n-1\}$, with $\zeta_i$ defined by (\ref{eq:zeta}).\hfill$\square$\vspace*{10pt}

\begin{remark}
Theorem~\ref{th:theo2} corrects a minor error in \cite{slotine}. Slotine proposed that the bounds for 
$\tilde{x}^{(i)}$ could be summarized as $|\tilde{x}^{(i)}|\le2^i\lambda^{i-n+1}\phi,i=0,1,\ldots,n-1$. 
Although both results lead to same bounds for $\tilde{x}$ and $\dot{\tilde{x}}$, they start to differ 
from each other when the order of the derivative is higher than one, $i>1$. For example, according to 
Slotine the bounds for the second derivative would be $|\ddot{\tilde{x}}|\le4\phi\lambda^{3-n}$ and 
not $|\ddot{\tilde{x}}|\le6\phi\lambda^{3-n}$, as demonstrated in Theorem~\ref{th:theo2}.
\end{remark}

\section{Concluding remarks}

In this work, a convergence analysis of smooth sliding mode controllers was presented. The 
attractiveness and invariant properties of the boundary layer as well as the exponential convergence  
of the tracking error vector to a bounded region were analytically proven. This last result corrected 
flawed conclusions previously reached in the literature. 

\section*{Acknowledgements}

The author acknowledges the support of the State of Rio de Janeiro Research Foundation (FAPERJ).
Furthermore, the author would like to thank Prof.\ Roberto Barrêto and Prof.\ Gilberto Corrêa for 
their insightful comments and suggestions.

\end{document}